\documentclass[11pt]{amsart}

\usepackage{amsmath}
\usepackage{amssymb}
\usepackage{bbm}
\usepackage{pdfsync}
\usepackage{esint} 
\usepackage[breaklinks=true, pagebackref]{hyperref}
\usepackage[utf8]{inputenc}
\usepackage[T1]{fontenc}
\usepackage{mathtools}
\usepackage{comment}

\renewcommand*{\backref}[1]{}
\renewcommand*{\backrefalt}[4]{%
	\ifcase #1 (Not cited.)%
	\or        (Cited on page~#2)%
	\else      (Cited on pages~#2)%
	\fi}

\usepackage[capitalize]{cleveref} 
\crefname{equation}{}{} 
\crefname{claim}{claim}{claims} 
\crefname{rem}{remark}{remarks} 

\crefname{page}{p.}{pp.}

\usepackage{enumitem}

\hypersetup{ 
	linktoc=page, 
	colorlinks=true,
	linkcolor={red!90!black}, 
	citecolor={blue!90!black}, 
	urlcolor=magenta,
	filecolor=teal,
}

\let\oldtocsection=\tocsection
\let\oldtocsubsection=\tocsubsection
\let\oldtocsubsubsection=\tocsubsubsection
\renewcommand{\tocsection}[2]{\hspace{0em}\oldtocsection{#1}{#2}}
\renewcommand{\tocsubsection}[2]{\hspace{1em}\oldtocsubsection{#1}{#2}}
\renewcommand{\tocsubsubsection}[2]{\hspace{2em}\oldtocsubsubsection{#1}{#2}}

\makeatletter
\newcommand\@dotsep{4.5}
\def\@tocline#1#2#3#4#5#6#7{\relax
	\ifnum #1>\c@tocdepth 
	\else
	\par \addpenalty\@secpenalty\addvspace{#2}%
	\begingroup \hyphenpenalty\@M
	\@ifempty{#4}{%
		\@tempdima\csname r@tocindent\number#1\endcsname\relax
	}{%
		\@tempdima#4\relax
	}%
	\parindent\z@ \leftskip#3\relax
	\advance\leftskip\@tempdima\relax
	\rightskip\@pnumwidth plus1em \parfillskip-\@pnumwidth
	#5\leavevmode\hskip-\@tempdima #6\relax
	\leaders\hbox{$\m@th
		\mkern \@dotsep mu\hbox{.}\mkern \@dotsep mu$}\hfill
	\hbox to\@pnumwidth{\@tocpagenum{#7}}\par
	\nobreak
	\endgroup
	\fi}
\makeatother 

\newcommand{\M}{{\mathcal M}}       %
\newcommand{\R}{{\mathbb R}}       

\newcommand{\HH}{{\mathcal H}}

\newcommand{\rf}[1]{{(\ref{#1})}}

\newcommand{\supp}{\operatorname{supp}}

\usepackage{dsfont}

\DeclareMathOperator{\Span}{span}
\NewDocumentCommand{\sff}{}{\mathrm{I\!I}}

\newcommand{\rom}[1]{%
	\textup{\uppercase\expandafter{\romannumeral#1}}%
}

\def\XXint#1#2#3{{\setbox0=\hbox{$#1{#2#3}{\int}$ }
		\vcenter{\hbox{$#2#3$ }}\kern-.58\wd0}}

\makeatletter
\newcommand{\doublewidetilde}[1]{{%
		\mathpalette\double@widetilde{#1}%
}}
\newcommand{\double@widetilde}[2]{%
	\sbox\z@{$\m@th#1\widetilde{#2}$}%
	\ht\z@=.8\ht\z@
	\widetilde{\box\z@}%
}
\makeatother

\usepackage{todonotes}

\newtheorem{theorem}{Theorem}[section]
\newtheorem{lemma}[theorem]{Lemma}

\newtheorem*{claim*}{Claim}
\newtheorem{conjecture}{Conjecture}
\newtheorem*{theorem*}{Theorem}

\newtheorem*{coro*}{Corollary}

\theoremstyle{definition}

\numberwithin{equation}{section}

\allowdisplaybreaks

\definecolor{taronja}{rgb}{0.9,0.5,0.05}
\usepackage{xcolor}

\begin{document}
\title[Unique continuation for locally uniformly distributed measures]{Unique continuation for locally uniformly distributed measures}

\author[Max Engelstein]{Max Engelstein}
\address{Max Engelstein, School of Mathematics, University of Minnesota, Minneapolis, MN, 55455, USA.}
\email{mengelst@umn.edu}

\author[Ignasi Guillén-Mola]{Ignasi Guillén-Mola}
\address{Ignasi Guillén-Mola, Departament de Matem\`atiques, Universitat Aut\`onoma de Barcelona.}
\email{ignasi.guillen@uab.cat}

\date{\today} 


\thanks{\textit{Thanks}: The authors would like to thank the referee for a careful reading and helpful comments. M.E. would like to thank Carlos Kenig for introducing him to \cref{con:locallyanalytic} many years ago. I.G.M. would like to thank the first author, Max Engelstein, for his hospitality and mentorship during I.G.M.'s stay at the University of Minnesota (UMN). This work began when I.G.M. was visiting M.E. at the UMN, during Spring 2024. Both authors would like to thank the university for its hospitality. M.E. was partially supported by NSF DMS CAREER 2143719, and the George W. Taylor Career Development Award. I.G.M. was supported by Generalitat de Catalunya’s Agency for Management of University and Research Grants (AGAUR) (2021 FI\_B 00637 and 2023 FI-3 00151), and also partially supported by 2021-SGR-00071 (AGAUR, Catalonia) and by the Spanish State Research Agency (AEI) projects PID2020-114167GB-I00 and PID2021-125021NAI00. I.G.M. thanks CERCA Programme/Generalitat de Catalunya for institutional support.
}


\begin{abstract}
    In this note we show that the support of a locally $k$-uniform measure in $\mathbb R^{n+1}$ satisfies a kind of unique continuation property. As a consequence, we show that locally uniformly distributed measures satisfy a weaker unique continuation property. This continues work of Kirchheim and Preiss in \cite{Kirchheim-Preiss-2002} and David, Kenig and Toro in \cite{David-Kenig-Toro-2001} and lends additional evidence to the conjecture proposed by Kowalski and Preiss in \cite{Kowalski-Preiss-1987} that each connected component of the support of a locally $n$-uniform measure in $\mathbb R^{n+1}$ is contained in the zero set of a quadratic polynomial.
\end{abstract}

\maketitle



\section{Introduction}\label{sec:intro}

In this note we consider locally uniform measures, first introduced by Kowalski and Preiss in \cite{Kowalski-Preiss-1987} and later studied by David, Kenig and Toro in \cite{David-Kenig-Toro-2001}. For $0\leq k \leq n+1$ integer, a Radon measure $\mu$ with $\supp \mu \subset \R^{n+1}$ is {\bf locally $k$-uniform} if \begin{equation}\label{localuniformcondition}
    \mu(B(x,r)) = w_k r^k \text{ whenever } x\in\supp\mu\text{ and } 0<r\leq 1,
\end{equation}
where $w_k = \HH^k (B^k (0,1))$ is the measure of the unit ball $B^k (0,1)$ in $\R^{k}$. When $\supp \mu$ is a smooth manifold and $\mu = \HH^k|_{\supp \mu}$ it is a natural question in differential geometry to ask what the condition \rf{localuniformcondition} tell us about the manifold $\supp \mu$, see e.g., \cite{Karp-Pinsky-1989,Osserman-1975}. 

For general $\mu$, the condition \rf{localuniformcondition} and its global cousin (i.e., \rf{localuniformcondition} without the restriction that $r \leq 1$) are extremely well studied, as they represent ``end point cases'' in the study of measure densities, see e.g., \cite{Preiss-1987,Marstrand-1964}. These measure density theorems have also been applied to understand rectifiability in measure spaces (see \cite{Camillo-2008-book,Merlo-2022}) and higher notions of smoothness, see \cite{David-Kenig-Toro-2001,Preiss-Tolsa-Toro-2009,Tolsa-2015-JLMS}.

The condition \rf{localuniformcondition}, is extremely rigid and constrains not just the measure, but also its support. For example any two locally uniform measures (in fact locally uniformly distributed, see \cref{localunifdistrcondition}) with the same support must be scalar multiples of one another, see for instance \cite[Theorem 3.4]{Mattila1995}\footnote{The theorem is written for uniformly distributed measures but the proof only uses the local structure of the measures.}. Furthermore, by Marstrand's theorem, \cite{Marstrand-1964}, there are no locally $k$-uniform measures for non-integer $k$.

For globally uniform measures (\rf{localuniformcondition} without the restriction $r\leq 1$), even more is known. It was proven by Kowalski and Preiss in \cite{Kowalski-Preiss-1987} that globally $k$-uniform measures in $\R^{n+1}$ must be of the form $\HH^k|_{\Sigma}$ where $\Sigma$ is a quadratic variety. To be more precise, in \cite[Theorem 4.1]{Kowalski-Preiss-1987}, it is only proven that in codimension $1$ the support is a subset of a quadratic variety. The proof for arbitrary codimension works mutatis mutandis, see \cite[Remark 17.4]{Mattila1995}, and can be found in \cite[Theorem 17.3]{Mattila1995}. That the measure is equal to $\HH^k|_\Sigma$ then follows from \cite[Theorem 1.4]{Kirchheim-Preiss-2002}. When $k = n$ (the codimension $1$ case) the admissible quadratic varieties were fully characterized by Kowalski and Preiss in \cite{Kowalski-Preiss-1987}, with $\Sigma$ given by (after a translation and rotation)
$$
\begin{aligned}
\Sigma &= \{(x_1,\ldots,x_{n+1})\in\R^{n+1} : x_{n+1}=0\},\text{ or}\\
\Sigma &= \{(x_1,\ldots,x_{n+1})\in\R^{n+1} : x_4^2=x_1^2+x_2^2+x_3^2\} \text{ if }n\geq 3.
\end{aligned}
$$
In higher codimension classifying globally $0$-uniform and globally $1$-uniform measures is straightforward and it follows from work in \cite{Kowalski-Preiss-1987} that any connected component of the support of a globally $2$-uniform measure must be contained in a sphere or a plane\footnote{According to \cite{David-Kenig-Toro-2001}, in unpublished work Kowalski showed that the smooth support of any locally $2$-uniform measure must be a union of spheres and planes in $\R^3$, see also \cite[Theorem 3.2(b)]{Kowalski-Preiss-1987} and \cite[Proposition 3.1]{Karp-Pinsky-1989}. Also note, that while a single sphere is not the support of a globally $2$-uniform measure, it is hard to rule out the possibility that countably many disjoint spheres, carefully placed, support a globally $2$-uniform measure, see \cite{Nimer-2022-JDG} for interesting computations in this direction.}. However, the characterization of which quadratic varieties support a globally $k$-uniform measures when $2<k<n$ remains an open problem and even constructing interesting examples is extremely difficult (see \cite{Nimer-2022-JDG}). 

However, given their success in classifying globally $n$-uniform measures in $\R^{n+1}$ it was natural for Kowalski and Preiss, \cite[Conjecture 5.14]{Kowalski-Preiss-1987}, to conjecture that each connected smooth hypersurface $\M\subset\R^{n+1}$ satisfying $\HH^n (\M\cap B(x,r))=w_n r^n$ for each $x\in \M$ and each sufficiently small $r>0$ is a subset of an $n$-dimensional quadric. We reformulate this conjecture in terms of locally $n$-uniform measures in $\R^{n+1}$:

\begin{conjecture}[See {\cite[Conjecture 5.14]{Kowalski-Preiss-1987}}]\label{con:locallyanalytic}
     Let $\mu$ be a locally $n$-uniform measure in $\R^{n+1}$. Then each connected component of the support of $\mu$ is contained in the zero set of a quadratic polynomial.
\end{conjecture}

This conjecture remains completely open for $n\geq 3$ with little progress (apart from work of \cite{David-Kenig-Toro-2001} which we will discuss below). In this note, provide some evidence towards \cref{con:locallyanalytic}, \emph{in any codimension}. In particular, we prove that the support of a locally $k$-uniform measure in $\R^{n+1}$ must satisfy the following ``strong unique continuation property'' in terms of the mean curvature vector:

\begin{theorem}[Strong Unique Continuation]\label{t:SUCP}
    Let $\mu$ be a locally $k$-uniform measure in $\R^{n+1}$. Assume there is a point $z_0 \in \supp\mu$ such that the mean curvature vector of $\supp\mu$ vanishes at $z_0$. Then the connected component of $\supp \mu$ containing $z_0$ must be a $k$-plane.
\end{theorem}

We prove this theorem in \cref{s:UCP}, where we also recall the definition of the mean curvature vector and establish a seemingly novel and key identity, \cref{lemma:key equality}. The assumption that the mean curvature vector vanishes (and is, in particular, well-defined) requires $\supp\mu$ be of class $C^2$ in a neighborhood of $z_0\in\supp\mu$. For $k<n$, all that is known (to the best of the authors' knowledge) is that $\supp\mu$ is a $k$-dimensional $C^{1,\beta}$-manifold ($0<\beta<1$) away from a closed set $\mathcal S$ such that $\HH^k(\mathcal S)=0$, see \cite[Theorem 1.9]{Preiss-Tolsa-Toro-2009}. For $k=n$, work of \cite[Theorem 1.10]{David-Kenig-Toro-2001} by David, Kenig and Toro, combined with \cite[Corollary 1.11]{Preiss-Tolsa-Toro-2009} by Preiss, Tolsa and Toro, shows $\supp \mu = \mathcal R\cup \mathcal S$, where $\dim \mathcal S \leq n-3$ (discrete if $n=3$) and around every $z_0\in \mathcal R$ there exists a $r_0>0$ such that $B(z_0,r_0) \cap \supp \mu$ is the graph of a $C^\infty$ function. Thus in codimension $1$ asking that the mean curvature vector exists is not particularly restrictive, while in higher codimensions it may in fact place additional constraints on $\mu$ (though it seems more likely that a similar regularity theory is true in the higher codimension case). 

Observe that a standard density argument\footnote{See for instance the proof of \cite[Theorem 4.5]{Kowalski-Preiss-1987}, or the proof of \cref{c:WUCP}, more precisely \rf{mu is in fact hausdorff}, in \cref{s:WUCP}.} shows that, if $\mu$ is locally $k$-uniform, then in any ball, $B$, such that $\supp \mu\cap B$ is of class $C^1$, we have $\mu|_{B} = \HH^k|_{\supp \mu\cap B}$. Therefore, \rf{localuniformcondition} (and the work of \cite{Preiss-Tolsa-Toro-2009}) immediately implies that any (locally) $k$-uniform measure $\mu$ can be written $\mu = \HH^k|_{\supp \mu}$. As a consequence, it is not more restrictive to study manifolds, $\M$, such that $\HH^k(\M\cap B(x,r)) = w_k r^k$ for all $x\in \M$ and $r>0$ sufficiently small. Additionally, in the setting of \cref{t:SUCP} above, we can conclude that $\mu|_P = \HH^k|_P$, where $P$ is the $k$-plane coinciding with the connected component of $\supp\mu$ containing $z_0$.

By the curvature restrictions on the support of (locally smooth) locally $n$-uniform measures in $\R^{n+1}$, see \cite[Corollary 3.1]{Kowalski-Preiss-1987}, \cref{t:SUCP}  still holds if the mean curvature is replaced by scalar curvature in codimension $1$.

\Cref{t:SUCP} provides some evidence towards \cref{con:locallyanalytic} (in any codimension), since if $\Sigma = \{Q = 0\}$ for some quadratic polynomial $Q$, and the second fundamental form of $\Sigma$ vanishes at a point, then $\Sigma$ must be a plane. We should mention that we know of only one other result regarding the mean curvature of the support of a locally $n$-uniform measure in $\R^{n+1}$: \cite[Theorem 3.3]{Kowalski-Preiss-1987} shows, under the much stronger assumption that $\supp \mu$ is minimal (i.e., the mean curvature vanishes everywhere) that $\supp \mu$ must be a union of planes. 

As a simple corollary, \cref{t:SUCP} implies a ``weak'' unique continuation type property for the support of a broader class of measures, which we call locally uniformly distributed. A Radon measure $\mu$ with $\supp \mu \subset \R^{n+1}$ is {\bf locally uniformly distributed} if
\begin{equation}\label{localunifdistrcondition}
0<\mu (B(x,r)) = \mu(B(y,r)) <\infty \text{ whenever }x,y\in \supp\mu \text{ and } 0<r\leq 1.
\end{equation}

Every locally uniform measure is locally uniformly distributed but the opposite is not always the case. Indeed, $\HH^3|_{\mathbb S^3}$ is locally uniformly distributed but a straightforward computation shows that it is not locally $3$-uniform in $\R^4$ (one can notice that $\mathbb S^3$ does not satisfy the curvature conditions of \cite[Corollary 3.1]{Kowalski-Preiss-1987}). In contrast $\HH^2|_{\mathbb S^2}$ is locally $2$-uniform in $\R^3$ (and locally uniformly distributed). Locally uniformly distributed measures naturally arise as the conical part of globally uniform measures (\rf{localuniformcondition} without the restriction $r \leq 1$), and so understanding their structure is an important step towards producing more examples of globally uniform measures, see e.g., \cite{Nimer-2022-JDG}.

Very little is known even about globally uniformly distributed measures (\rf{localunifdistrcondition} without the restriction that $r \leq 1$). Kirchheim and Preiss showed in \cite{Kirchheim-Preiss-2002} that globally uniformly distributed measures are supported on analytic varieties. Given the few known examples, it is natural to conjecture that the support of a locally uniformly distributed measure must also be an analytic variety. Our next theorem supports this conjecture.

\begin{theorem}[Weak Unique Continuation]\label{c:WUCP}
    Let $\mu$ be a locally uniformly distributed measure in $\R^{n+1}$. Assume there is an open ball $B(z_0,r_0)$ (with $z_0\in \supp\mu$ and $r_0>0$) and a $k$-plane $P$ such that $B(z_0,r_0)\cap \supp\mu = B(z_0,r_0)\cap P$. Then there is a unique connected component of $\supp \mu$ which intersects non-trivially with $B(z_0, r_0)$ and this connected component is the $k$-plane $P$.
\end{theorem}

This follows from the strong unique continuation \cref{t:SUCP} for locally $k$-uniform measures in $\R^{n+1}$, and its proof is given in \cref{s:WUCP}.

\section{Proof of \texorpdfstring{\cref{c:WUCP}:}{} Weak unique continuation property}\label{s:WUCP}

Before proving \cref{t:SUCP}, we show how the weak unique continuation property in \cref{c:WUCP} follows from our \cref{t:SUCP}.

\begin{proof}[Proof of \cref{c:WUCP}]
    We mainly follow the first part of \cite[Proof of Theorem 1.4(ii)]{Kirchheim-Preiss-2002}. The Lebesgue differentiation theorem implies that the limit
    $$
    K(y) \coloneqq \lim_{r\to 0} \frac{\mu(B(y,r))}{(\mu + \HH^k|_{\supp \mu})(B(y,r))}
    $$
    exists at $(\mu + \HH^k|_{\supp \mu})$-a.e.\ $y\in\supp\mu$, and that
    \begin{equation}\label{eq:mu is integral limit of mu and Lebesgue}
    \mu (E) = \int_{E} K(y) \, d(\mu + \HH^k|_{\supp \mu})(y) \text{ for any } E\subset B(z_0,r_0),
    \end{equation}
    since $\mu$ is (clearly) absolutely continuous with respect to $\mu + \HH^k|_{\supp \mu}$, i.e., $\mu \ll (\mu + \HH^k|_{\supp \mu})$. Note also that the lower bound in \rf{localunifdistrcondition} implies $K(y)>0$ for $(\mu + \HH^k|_{\supp \mu})$-a.e.\ $y\in\supp\mu$.

    Fixed a (any) point $p\in \supp\mu$, we define
    $$
    g_\mu (r) \coloneqq \mu(B(p,r))\text{ for }0<r\leq 1,
    $$
    see \rf{localunifdistrcondition}. By the flatness assumption $B(z_0,r_0)\cap\supp\mu = B(z_0,r_0)\cap P$ with $k$-plane $P$, for $y\in B(z_0,r_0)\cap \supp\mu$ we have
    $$
    K(y) = \lim_{r\to 0} \frac{g_\mu(r)}{g_\mu(r) + w_k r^k} \leq 1,
    $$
    which is independent of $y\in B(z_0,r_0)\cap\supp\mu$. Since $K(y)$ exists $(\mu + \HH^k|_{\supp \mu})$-a.e., we conclude that the limit $K(y)$ exists for all $y\in B(z_0,r_0)\cap\supp\mu$. Furthermore, since $K(y) > 0$ for $(\mu + \HH^k|_{\supp \mu})$-a.e.\ $y\in\supp\mu$, it follows that $0 < K(y) \leq 1$ for all $y\in B(z_0,r_0)\cap\supp\mu$.
    
    We define (before checking it makes sense)
    $$
    c \coloneqq \lim_{r\to 0}\frac{g_\mu (r)}{w_k r^k}.
    $$
    We have that this limit exists by the existence of the limit $K(y)$ for all $y\in B(z_0,r_0)\cap\supp\mu$ and since
    $$
    \frac{1}{c} = \lim_{r\to 0}\frac{w_k r^k}{g_\mu(r)} = \lim_{r\to 0} \frac{w_k r^k+g_\mu(r)}{g_\mu(r)} - 1 = \frac{1}{K(y)} - 1.
    $$
    Moreover, the fact that $K(y) > 0$ implies that $c > 0$. Furthermore if $c = \infty$ then $K(y) = 1$ for $y\in B(z_0, r_0) \cap \supp \mu$ and therefore \rf{eq:mu is integral limit of mu and Lebesgue} would give $\HH^k|_{\supp\mu}(E)=0$ for any $E\subset B(z_0,r_0)\cap \supp \mu$, but this cannot happen since $B(z_0,r_0)\cap\supp\mu = B(z_0,r_0)\cap P$.
    
    Since $c\in (0,\infty)$, we can write $K(y) = \frac{c}{c+1}\in (0,1)$ for all $y\in B(z_0,r_0)\cap\supp\mu$, and therefore by \rf{eq:mu is integral limit of mu and Lebesgue} we conclude that
    \begin{equation}\label{mu is in fact hausdorff}
    \mu (E) = c \HH^k|_{\supp \mu}(E) \text{ for any } E\subset B(z_0,r_0).
    \end{equation}
    
    The equality \rf{mu is in fact hausdorff} in particular implies $g_\mu(r) = \mu(B(x,r)) = cw_k r^k$ for all $x\in B(z_0,r_0/2)\cap\supp\mu$ and $0<r\leq r_0/4$. However, by the locally uniformly distributed condition \rf{localunifdistrcondition}, the same holds for any $x\in \supp\mu$. That is,
    $$
    g_\mu(r) = \mu(B(x,r)) = cw_k r^k \text{ for all }x\in \supp\mu \text{ and } 0<r\leq r_0/4.
    $$
    Thus, the measure
    $$
    \widetilde \mu (\cdot) \coloneqq \frac{1}{c}\frac{1}{\left(\frac{r_0}{4}\right)^k}g\left(\frac{r_0}{4} \cdot\right)
    $$
    satisfies the local uniform condition \rf{localuniformcondition}, is flat in the open ball $B(4z_0/r_0,4)=B(z_0,r_0)/(r_0/4)$, and $\supp \widetilde{\mu}$ is a dilation of $\supp \mu$. The result follows by \cref{t:SUCP} applied to the locally $k$-uniform measure $\widetilde\mu$.
\end{proof}

\section{Proof of \texorpdfstring{\cref{t:SUCP}:}{} Strong unique continuation property}\label{s:UCP}

In this section we prove the strong continuation-type property in \cref{t:SUCP} for locally $k$-uniform measures in $\R^{n+1}$ that is, if the mean curvature vector of the support vanishes at a point, then the whole connected component containing the point must be a $k$-plane.

Let us recall that the mean curvature vector of a $k$-dimensional $C^2$-manifold in $\R^{n+1}$ is the trace of the (vector-valued) second fundamental form, divided by $k$. Let $\M\cap B(z_0,r_0)$, with $z_0\in \M$ and $r_0>0$, be a $C^2$-manifold of dimension $k$. If $\{\hat{n}_j (z)\}_{j=k+1}^{n+1}$ is an orthonormal basis of the normal space $N_\M(z)$ at $z\in\M\cap B(z_0,r_0)$ then, intuitively, the $j$-component of the mean curvature vector is a measure of ``flatness'' in that it infinitesimally captures how much the measure of $\M$ changes under perturbations in the direction of $\hat{n}_j(z)$.

For our purposes it is enough to work in coordinates, where the definition becomes more concrete. Assume that $\M\cap B(z_0,r_0)$ can be described by a $C^2$ function $\varphi=(\varphi_{k+1},\ldots,\varphi_{n+1}):\R^k\to \R^{n+1-k}$ such that $z=(z^\prime,\varphi(z^\prime))$ for $z\in \M\cap B(z_0,r_0)$, with $z^\prime \in \R^k$. Note that, in this case, the tangent space $T_\M (z)$ is
$$
T_\M (z) \coloneqq \Span \left(\left\{ t_j \right\}_{j=1}^k\right), \text{ with } t_j(z) \coloneqq \frac{e_j^k+\partial_{x_j} \varphi(z^\prime)}{\sqrt{1+|\partial_{x_j} \varphi(z^\prime)|^2}},
$$
and $\{e_j^k\}_{j=1}^k$ is the standard orthonormal basis of $\R^k$. As above, $\{\hat{n}_j(z)\}_{j=k+1}^{n+1}$ is an orthonormal basis of the normal space $N_\M(z)$. Then the vector-valued second fundamental form $\sff_{\M}(\cdot,\cdot)(z):T_\M(z)\times T_\M (z)\to N_\M (z)$ is given by
$$
\sff_{\M} (u,v)(z) \coloneqq \Pi_{N_\M (z)}(u\cdot \nabla v) = \sum_{j=k+1}^{n+1}\left\langle u\cdot \nabla v, \hat{n}_j(z)\right\rangle \hat{n}_j(z).
$$
(Where the inner product is the Euclidean inner product).

The mean curvature vector is then the trace so we have that
$$
H_\M(z) = \frac{1}{k} \sum_{i,j=1}^k g^{ij}(z)\sum_{\ell = k+1}^{n+1} \left\langle t_j(z) \cdot \nabla t_i(z), \hat{n}_\ell(z)\right\rangle \hat{n}_\ell(z),
$$
where $g^{ij}(z)$ is the inverse of the first fundamental form, given by $g_{ij}(z) = \left\langle t_i(z), t_j(z)\right\rangle$. 

In generality, this expression is extremely complicated, but it simplifies greatly if we assume that $T_{\M}(0) = \R^k\times\{0\}^{n+1-k}$ (that is, $\partial_{x_i} \varphi(0) = 0$, for every $i=1,\ldots, k$). In this setting we can take $t_j(0) = e_j$, $j=1,\ldots,k$, and $\hat{n}_j(0) = e_j$, $j=k+1,\ldots,n+1$ (here $\{e_i\}_{i=1}^{n+1}$ is the standard orthonormal basis of $\R^{n+1}$). It follows that $g^{ij}(0) = \delta_{ij}$ (i.e., is 1 when $i= j$ and $0$ otherwise). Thus our formula simplifies to $$H_\M(0) = \frac{1}{k}\sum_{i=1}^k \sum_{\ell = k+1}^{n+1} \left\langle e_i \cdot \nabla t_i(0), e_\ell\right\rangle e_\ell.$$ We recall that $\nabla t_i(0)$ is a matrix and we need only compute its $(\ell, i)$-entry. This is a straightforward computation that gives
$$
(\nabla t_i(0))_{\ell, i} = \partial_{x_i} \left(\frac{\partial_{x_i} \varphi_\ell}{\sqrt{1+|\partial_{x_i} \varphi|^2}}\right)(0)\\
= \partial^2_{x_i}\varphi_\ell(0),
$$
where we have used in the second equality that $\partial_{x_i} \varphi(0) = 0$ for all $i= 1\ldots, k$.  Thus we have the final form
\begin{equation}\label{e:MCflat}
    H_\M(0) 
    = \frac{1}{k}\sum_{\ell=k+1}^{n+1}\Delta \varphi_\ell(0)e_\ell
    = \frac{1}{k} (0,\Delta\varphi(0)), \text{ if } T_\M(0)= \R^k\times\{0\}^{n+1-k}.
\end{equation}

Before we begin our proofs, let us introduce some preliminary notation. We define the vector $b_{z,r}$ and its nonnormalized version $\widetilde b_{z,r}$ as
\begin{subequations}
    \begin{align}
        b_{z,r} &\coloneqq \frac{k+2}{2w_k r^{k+2}} \int_{B(z,r)} \left(r^2-|y-z|^2\right)(y-z)\, d\mu(y),\label{vector b}\\
        \widetilde b_{z,r} &\coloneqq \frac{2w_k r^{k+2}}{k+2} b_{z,r}
        =\int_{B\left(z, r\right)}\left(r^2-|y-z|^2\right)\left(y-z\right) \, d \mu(y),\label{vector b normalized}
    \end{align}
\end{subequations}
and the quadratic form $Q_{z,r}$ and its normalized version $\widetilde Q_{z,r}$ on $\R^{n+1}$ as
\begin{subequations}
    \begin{align}
        Q_{z,r} (x) &\coloneqq \frac{k+2}{w_k r^{k+2}} \int_{B(z,r)} \langle x,y-z \rangle^2\, d\mu(y),\label{quadratic form Q}\\
        \widetilde  Q_{z,r} (x) &\coloneqq \frac{w_k r^{k+2}}{k+2} Q_{z,r}
        = \int_{B(z,r)} \langle x,y-z \rangle^2\, d\mu(y)\label{quadratic form Q normalized},
    \end{align}
\end{subequations}
all first introduced in \cite[(4.3) and (4.4)]{Kowalski-Preiss-1987}, see also \cite[(7.2) and (7.3)]{David-Kenig-Toro-2001}. Throughout this section, $b_{z,r}^j$ and $\widetilde b_{z,r}^j$ denote the $j$-component ($1\leq j \leq n+1$) of the vectors $b_{z,r}$ and $\widetilde b_{z,r}$ respectively.

Here we present the key equality in the proof of our main result \cref{t:SUCP}. This relates the orthogonal (see \rf{statement orthogonal}) vector $\widetilde b_{0,r}$ and the mean curvature (given in local coordinates by $\Delta \varphi(0)$, see \rf{e:MCflat}), to the integral of ``height'' of the support squared. This identity is the main novel part of our work. 

\begin{lemma}[Key equality]\label{lemma:key equality}
    Let $\mu$ be a locally $k$-uniform measure in $\R^{n+1}$ and $0<r<1/2$. Assume $0\in \supp\mu$ and that $\supp\mu\cap B(0,r_0)$ (for some $r_0>0$) is given by the graph of a $C^2$ function $\varphi = (\varphi_{k+1},\ldots,\varphi_{n+1}) : \R^k \to \R^{n+1-k}$ with $|\nabla\varphi_j(0)|=0$ for all $j=k+1,\ldots,n+1$. Then there holds
    $$
    \frac{1}{2} \sum_{j=k+1}^{n+1} \widetilde b_{0,r}^j \Delta \varphi_j (0) = \sum_{j=k+1}^{n+1} \int_{B(0,r)} y_j^2\, d\mu(y).
    $$
\end{lemma}

We prove this in \cref{sec:proof 2 key equality}. We are now ready to turn to the proof of \cref{t:SUCP}.

\begin{proof}[Proof of \cref{t:SUCP}]
    Since we are $C^2$-smooth in a neighborhood of $z_0$, we can work in coordinates: after a harmless rotation and translation, let $\supp\mu$ in a neighborhood of $0\in\supp\mu$ be given by $(x^\prime, \varphi(x^\prime))$, $x^\prime \in \R^n$, with $\varphi(x)=(\varphi_{k+1}(x^\prime),\ldots,\varphi_{n+1}(x^\prime))$ satisfying $\varphi(0)=0$ and $|\nabla\varphi_j(0)| = 0$ for all $j=k+1,\ldots,n+1$. As we computed in \rf{e:MCflat} above, this implies that 
    $$
    H_{\supp\mu} (0) = \frac{1}{k} (0,\Delta\varphi(0)).
    $$
    
    Fixed $0<r<1/2$, by \cref{lemma:key equality} we have
    $$
    \frac{1}{2} \sum_{j=k+1}^{n+1} \widetilde b_{0,r}^j \Delta \varphi_j (0) = \sum_{j=k+1}^{n+1} \int_{B(0,r)} y_j^2\, d\mu(y).
    $$
    Since $(0,\Delta \varphi (0))=0$ (in particular $\Delta \varphi_j (0)=0$ for all $j=k+1,\ldots,n+1$), then the integral in the right-hand side is zero, implying that $y_j=0$ for $\mu$-a.e.\ $y\in B(0,r)$ for all $j=k+1,\ldots,n+1$. But this is only possible if $B(0,r)\cap\supp\mu = B(0,r)\cap P$ for the $k$-plane $P\coloneqq\R^k\times\{0\}^{n+1-k}$. Since this holds for every $0<r<1/2$, we conclude that $\supp\mu\cap B(0,1/2) = P\cap B(0,1/2)$. 
    
    To show that $B(z_0,1/2)\cap \supp \mu= B(z_0,1/2)\cap P$ for any $z_0$ in the same connected component of $\supp \mu$, consider a chain of balls of radius $1/2$ centered at points in $\supp \mu$ such that the center of each new ball is in the previous ball with $0$ as the center of the first ball and $z_0$ as the center of the last ball. The result follows by applying this argument to each ball in the chain in sequence.
\end{proof}

\subsection{Proof of Key equality\texorpdfstring{: \cref{lemma:key equality}}{}}\label{sec:proof 2 key equality}

Recall the definition of $b_{z,r}$ and $Q_{z,r}$ from \rf{vector b} and \rf{quadratic form Q}. Let us state two fundamental properties satisfied by a locally $k$-uniform measure $\mu$ in $\R^{n+1}$. The first is the elementary identity
\begin{equation}\label{uniform measures integral 0}
    \int_{B(z,r)} |y-z|^2\, d\mu(y) = \frac{k w_k}{k+2} r^{k+2} \text{ for all }z\in \supp\mu\text{ and }0<r<1,
\end{equation}
for the detailed straightforward computation (for $k=n$) see \cite[(4.5)]{Kowalski-Preiss-1987} for instance. Secondly, proved in \cite[(10.3)]{David-Kenig-Toro-2001}\footnote{The codimension $1$ version is proved in \cite[(10.3)]{David-Kenig-Toro-2001}, but its proof does not depend on the codimension. For a statement in any codimension, see \cite[Proposition 2.3 (2.12)]{Preiss-Tolsa-Toro-2009}.}, for $z\in\supp\mu$ and $0<r<1/2$ there holds
\begin{equation}\label{inequality b Q and modulus}
\left|\langle 2b_{z,r},x-z\rangle + Q_{z,r}(x-z) - |x-z|^2\right| \leq C \frac{|x-z|^3}{r} \text{ for all } x\in \supp\mu \cap B(z,r/2).
\end{equation}
Since $Q_{z,r}(x-z)\leq \frac{k+2}{w_k}|x-z|^2$, this in particular implies that for $|x-z|$ small enough there holds
$$
\left|\langle 2b_{z,r},x-z\rangle\right| \leq C^\prime |x-z|^2.
$$
From this latter inequality we conclude that, if $\supp\mu\cap B(z,r_0)$ (for some $r_0>0$) is at least of class $C^1$, then either $b_{z,r}=0$ or
\begin{equation}\label{statement orthogonal}
    b_{z,r}\text{ is orthogonal to the tangent }k\text{-plane of }\supp\mu\text{ at }z\in \supp\mu.
\end{equation}

Having recalled these facts we can now establish the following algebraic relation, which gets us most of the way to \cref{lemma:key equality}. We thank the reviewer for pointing out that similar estimates (in spirit) actually date back at least to the work of \cite{Marstrand-1964}, see Lemma 4 there. We should also state explicitly that the general idea of working in local coordinates and estimating errors using the smoothness of the support is ubiquitous in the literature, appearing in \cite{Kowalski-Preiss-1987, Karp-Pinsky-1989, Mattila1995, David-Kenig-Toro-2001, Preiss-Tolsa-Toro-2009} and surely elsewhere.

\begin{lemma}\label{lemma:equality for every e in sphere}
    Under the assumptions of \cref{lemma:key equality}, there holds
    \begin{equation}\label{eq:equality for every e in sphere}
    \sum_{j=k+1}^{n+1} b_{0,r}^j \langle\nabla^2 \varphi_j(0) e, e\rangle = 1 - Q_{0,r}(e) \text{ for all } e\in \mathbb S^{k-1}\times \{0\}^{n+1-k}.
    \end{equation}
    \begin{proof}
        During the proof, by the $C^2$ local representation of $\supp\mu$ we write $x=(x^\prime, \varphi(x^\prime))$ with $x^\prime \in \R^k$ and $\varphi(x^\prime) = (\varphi_{k+1}(x^\prime),\ldots,\varphi_{n+1}(x^\prime)) \in \R^{n+1-k}$. As we are assuming that $0\in \supp\mu$ and $|\nabla\varphi_j(0)|=0$ for all $j=k+1,\ldots,n+1$, by Taylor's theorem we then have
        $$
        \varphi_j(x^\prime) = \frac{\langle\nabla^2 \varphi_j(0) x^\prime, x^\prime\rangle}{2} + o(|x^\prime|^2) \text{ for all } j=k+1,\ldots,n+1.
        $$
        (We use little $o$ notation, so that $o(s)$ is a quantity that when divided by $s$ goes to zero as $s\downarrow 0$).
        
        With $z=0$, \rf{inequality b Q and modulus} reads as
        \begin{equation}\label{inequality b Q and modulus restatement}
        \left|\langle 2b_{0,r},x\rangle + Q_{0,r}(x) - |x|^2\right| \leq C \frac{|x|^3}{r} \text{ for all } x\in \supp\mu \cap B(0,r/2).
        \end{equation}
        First, from \rf{statement orthogonal} we have that $b_{0,r}$ is orthogonal to the tangent $k$-plane of $\supp\mu$ at $0$, or it may be zero. That is, in local coordinates, this means $b_{0,r} \in (\R^k\times\{0\}^{n+1-k})^\perp$, since $|\nabla \varphi_j (0)|=0$ for all $j=k+1,\ldots,n+1$. This implies $b_{0,r}^j=0$ for all $j=1,\ldots,k$, and therefore,
        $$
        \langle 2b_{0,r},x\rangle = 2 \sum_{j=k+1}^{n+1} b_{0,r}^j \varphi_j(x^\prime) 
        = \sum_{j=k+1}^{n+1} b_{0,r}^j \langle\nabla^2 \varphi_j(0) x^\prime, x^\prime\rangle + o (|x^\prime|^2).
        $$
        Second,
        $$
        \begin{aligned}
        Q_{0,r} (x) 
        &=\frac{k+2}{w_k r^{k+2}} \int_{B(0,r)} \langle x,y \rangle^2\, d\mu(y)\\
        &=\frac{k+2}{w_k r^{k+2}} \int_{B(0,r)} \left(\langle x^\prime,y^\prime \rangle + \sum_{j=k+1}^{n+1} \varphi_j(x^\prime)y_j\right)^2\, d\mu(y)\\
        &= \frac{k+2}{w_k r^{k+2}} \int_{B(0,r)} \langle x^\prime, y^\prime \rangle^2\, d\mu(y) + o(|x^\prime|^2).
        \end{aligned}
        $$
        Finally, simply $|x|^2 = |x^\prime|^2 + o(|x^\prime|^3)$. Also from this, we have that $|x|\approx |x^\prime|$ for $x\in \supp\mu\cap B(0,r/2)$ with small enough modulus.
        
        From \rf{inequality b Q and modulus restatement} and these estimates, there exists a $\theta\in (0,1)$ (depending on $\|\varphi\|_{C^2}$) such that  for all $x\in \supp\mu \cap B(0,r/2)$ with $|x| \leq \theta r$ we get
        \begin{equation}\label{e:almostontangentplane}
        \left| \sum_{j=k+1}^{n+1} b_{0,r}^j \langle\nabla^2 \varphi_j(0) x^\prime, x^\prime\rangle + \frac{k+2}{w_k r^{k+2}}\int_{B(0,r)} \langle x^\prime, y^\prime \rangle^2\, d\mu(y) - |x^\prime|^2\right| \leq \frac{1}{r}o(|x^\prime|^2).
        \end{equation}
        Note that we can simply write the middle term as $Q_{0,r} ((x^\prime,0))$.
        
        Given $e\in \mathbb S^{k-1}\times\{0\}^{n+1-k}$, let $\theta r/2> \varepsilon > 0$ and $x^\prime = \varepsilon e$ in \rf{e:almostontangentplane}. Dividing both sides of \rf{e:almostontangentplane} by $|x^\prime|^2 = \epsilon^2$ and letting $\varepsilon \downarrow 0$, we get \rf{eq:equality for every e in sphere}.
    \end{proof}
\end{lemma}

We now turn to the proof of \cref{lemma:key equality}.

\begin{proof}[Proof of \cref{lemma:key equality}]
    Multiplying the equality \rf{eq:equality for every e in sphere} in \cref{lemma:equality for every e in sphere} by $w_k r^{k+2}/(k+2)$ we have
    \begin{equation}\label{eq:multiple 1 of equality for every e in sphere}
    \frac{1}{2} \sum_{j=k+1}^{n+1} \widetilde b_{0,r}^j \langle\nabla^2 \varphi_j(0) e, e\rangle = \frac{w_k}{k+2} r^{k+2} - \widetilde Q_{0,r}(e) \text{ for all } e\in \mathbb S^{k-1}\times \{0\}^{n+1-k}.
    \end{equation}
    For $1\leq i \leq k$, taking $e=e_i\in \mathbb S^{k-1}\times \{0\}^{n+1-k}$ in \rf{eq:multiple 1 of equality for every e in sphere}, recall $\{e_i\}_{i=1}^{n+1}$ is the standard orthonormal basis of $\R^{n+1}$, we have
    \begin{equation}\label{eq:partial i i}
    \frac{1}{2} \sum_{j=k+1}^{n+1} \widetilde b_{0,r}^j \partial_{x_i,x_i} \varphi_j(0) = \frac{w_k}{k+2} r^{k+2} - \int_{B(0,r)} y_i^2 \, d\mu(y).
    \end{equation}
    Adding up \rf{eq:partial i i} for $1\leq i \leq k$ we get
    $$
    \frac{1}{2} \sum_{j=k+1}^{n+1} \widetilde b_{0,r}^j \Delta \varphi_j(0) = \frac{k w_k}{k+2} r^{k+2} - \sum_{i=1}^k \int_{B(0,r)} y_i^2 \, d\mu(y).
    $$
    The lemma follows from the above after applying the identity \rf{uniform measures integral 0}.
\end{proof}

\bibliographystyle{alpha}
\bibliography{references-phd.bib} 

\end{document}